\documentclass[reqno]{amsart}
\usepackage{hyperref}

\AtBeginDocument{{\noindent\small
\emph{Electronic Journal of Differential Equations},
Vol. 2013 (2013), No. 71, pp. 1--8.\newline
ISSN: 1072-6691. URL: http://ejde.math.txstate.edu or http://ejde.math.unt.edu
\newline ftp ejde.math.txstate.edu}
\thanks{\copyright 2013 Texas State University - San Marcos.}
\vspace{9mm}}

\begin{document}
\title[\hfilneg EJDE-2013/71\hfil Riesz bases]
{Riesz bases generated by the spectra of Sturm-Liouville problems}

\author[T. Harutyunyan, A. Pahlevanyan, A. Srapionyan \hfil EJDE-2013/71\hfilneg]
{Tigran Harutyunyan, Avetik Pahlevanyan, Anna Srapionyan}  % in alphabetical order

\address{Tigran Harutyunyan \newline
Faculty of Mathematics and Mechanics, Yerevan State University,
 1 Alex Manoogian, 0025, Yerevan, Armenia}
\email{hartigr@yahoo.co.uk}

\address{Avetik Pahlevanyan \newline
Faculty of Mathematics and Mechanics, Yerevan State University,
1 Alex Manoogian, 0025, Yerevan, Armenia}
\email{avetikpahlevanyan@gmail.com}

\address{Anna Srapionyan \newline
Faculty of Mathematics and Mechanics, Yerevan State University,
1 Alex Manoogian, 0025, Yerevan, Armenia}
\email{srapionyan.anna@gmail.com}

\thanks{Submitted November 19, 2012. Published March 17, 2013.}

\begin{abstract}
 Let $\{\lambda _n^2\} _{n = 0}^\infty$ be the spectra of a Sturm-Liouville
 problem on $[0,\pi ]$. We investigate the question:
 Do the systems $\{ \cos(\lambda_nx)\} _{n = 0}^\infty$ or
 $\{ \sin(\lambda_n x)\} _{n = 0}^\infty$ form Riesz bases in
 ${L^2}[0,\pi ]$? The answer is almost always positive.
\end{abstract}

\maketitle
\numberwithin{equation}{section}
\newtheorem{theorem}{Theorem}[section]
\newtheorem{lemma}[theorem]{Lemma}
\newtheorem{definition}[theorem]{Definition}
\newtheorem{remark}[theorem]{Remark}
\allowdisplaybreaks

\section{Introduction and statement of results}
\label{sec1}

Let ${\mu _n} = \lambda _n^2(q, \alpha, \beta)$, $n=0,1,2,\dots$
be the eigenvalues of the Sturm-Liouville
boundary-value problem $L(q, \alpha, \beta)$
\begin{gather} \label{eq1.1}
- y'' + q(x) y = \mu y, \quad x \in (0,\pi), \; \mu  \in\mathbb{C},\\
 \label{eq1.2}
y(0)\cos \alpha + y'(0)\sin \alpha = 0, \quad \alpha \in (0,\pi],\\
 \label{eq1.3}
y(\pi)\cos \beta + y'(\pi)\sin \beta = 0, \quad \beta \in [0,\pi),
\end{gather}
where $q \in {L_R^1}[0,\pi ]$, that is $q$ is a real, summable on $[0,\pi]$
function. In the simplest case, when $q(x) = 0$ almost everywhere (a.e.)
 on $[0,\pi]$, eigenfunctions of the problem $L(0, \alpha, \beta)$,
 which satisfy the initial conditions $y(0) = \sin \alpha$,
 $y'(0) = - \cos \alpha$, have the form
\begin{equation} \label{eq1.4}
\varphi _n^0(x, \alpha, \beta) = \cos (\lambda _n(0, \alpha, \beta)x)
\sin \alpha
- \frac{{\sin ( {{\lambda _n}(0, \alpha, \beta)x} )}}{{{\lambda _n}
(0, \alpha, \beta)}}\cos \alpha, \quad n=0,1,2,\dots
\end{equation}
and form an orthogonal basis in ${L^2}[0, \pi]$.
Here rises a natural question: Do the systems of functions
 $\{ \cos({\lambda _n}(0, \alpha, \beta)x)\} _{n = 0}^\infty$ and
  $\{ \sin({\lambda _n}(0, \alpha, \beta)x)\} _{n = 0}^\infty$
separately form basis in ${L^2}[0,\pi ]$?
 Examples show, that the answer is not always positive and depends on
 $\alpha$ and $\beta$. When $\alpha = \beta = \frac{\pi}{2}$, then
${\lambda _n}(0, \frac{\pi }{2}, \frac{\pi }{2}) = n$, $n=0,1,2,\dots$
and the system
$\{ \cos ({\lambda _n}(0, \frac{\pi}{2}, \frac{\pi}{2})x)\} _{n = 0}^\infty
= \{ \cos(nx)\} _{n = 0}^\infty$ forms an orthogonal basis,
but the system $\{ \sin ({\lambda _n}
(0, \frac{\pi }{2}, \frac{\pi }{2})x)\} _{n = 0}^\infty  = \{ 0\}
 \cup \{ \sin(nx)\} _{n = 1}^\infty$
is not a basis because of the ``unnecessary'' member $\sin (0x) \equiv 0$.
However, throwing away this ``unnecessary'' member, we obtain an orthogonal
basis $\{ \sin(nx)\} _{n = 1}^\infty$. In the case of $\alpha = \pi$,
$\beta = 0$
(see below section~\ref{sec2}), ${\lambda _n}(0, \pi, 0) = n + 1$,
$n=0,1,2,\dots$ and the system
$\{ \sin( {\lambda _n}(0, \pi, 0)x)\} _{n = 0}^\infty
= \{ \sin( (n + 1)x)\} _{n = 0}^\infty$
forms an orthogonal basis, but the system
$\{ \cos( {\lambda _n}(0, \pi, 0)x)\} _{n = 0}^\infty
= \{ \cos( (n + 1)x)\} _{n = 0}^\infty $ is not complete
in ${L^2}[0, \pi]$, there is a lack of constant, but adding it, thus,
taking the system
$\{ 1\}  \cup \{ \cos ((n + 1)x)\} _{n = 0}^\infty
 = \{ \cos(nx)\} _{n = 0}^\infty $
we obtain a basis in ${L^2}[0, \pi]$. The question that we want to answer
in this paper is the following: Do the systems
 $\{ \cos(\lambda _n(q, \alpha, \beta)x)\} _{n = 0}^\infty $ and
$\{ \sin(\lambda _n(q,\alpha,\beta)x)\} _{n = 0}^\infty $ form Riesz bases in
${L^2}[0,\pi ]$?
The answer we formulate in  theorems \ref{thm1.1} and \ref{thm1.2} below.

\begin{theorem} \label{thm1.1}
The system of functions
$\{ \cos( {\lambda _n}(q, \alpha, \beta)x)\} _{n = 0}^\infty $
is a Riesz basis in ${L^2}[0, \pi]$ for each triple
$( {q, \alpha, \beta} ) \in L_R^1[ {0, \pi} ] \times ( {0,\pi } ]
\times [ {0,\pi } )$, except one case: when $\alpha  = \pi, ~\beta = 0$,
the system $\{ \cos( {\lambda _n}(q, \pi, 0)x)\} _{n = 0}^\infty $
is not a basis, but the system
$\{ f(x)\}  \cup \{ \cos( {\lambda _n}(q, \pi, 0)x)\} _{n = 0}^\infty $
is a Riesz basis in ${L^2}[0, \pi]$, if $f(x)=\cos(\lambda x)$,
where $\lambda ^2 \neq \lambda _n^2$ for every $n=0,1,2,\dots$.
\end{theorem}

\begin{theorem} \label{thm1.2}
\begin{itemize}
\item[1.]
Let $\alpha$, $\beta  \in (0, \pi)$. Then the systems
\begin{itemize}
\item[(a)]
$\{ \sin (\lambda_n x)\} _{n = 1}^\infty $, if there is no zeros among
${\lambda _n} = {\lambda _n}(q, \alpha, \beta )$, $n=0,1,2,\dots$
(i.e. in this case we ``throw away'' $\sin(\lambda_0 x)$),

\item[(b)]
$\{ \sin(\lambda_n x)\} _{n = 0}^{{n_0} - 1} \cup
\{ \sin(\lambda_n x)\} _{n = {n_0} + 1}^\infty $, if
${\lambda _{n_0}}(q, \alpha, \beta ) = 0$ (we "throw away"
$\sin(\lambda_{n_0}x) \equiv 0$).
\end{itemize}
are Riesz bases in ${L^2}[0, \pi]$.
\item[2.]
Let $\alpha  = \pi$,  $\beta \in (0, \pi)$ or $\alpha \in (0, \pi)$,
$\beta = 0$. Then the systems
\begin{itemize}
\item[(a)]
$\{\sin(\lambda_n x)\} _{n = 0}^\infty$, if there is no zeros among
${\lambda _n} = {\lambda _n}(q, \alpha,\beta)$, $n=0,1,2,\dots$,
\item[(b)]
$\{ \sin(\lambda_n x)\} _{n = 0}^{{n_0} - 1} \cup \{ x\}
\cup \{ \sin(\lambda_n x)\} _{n = {n_0} + 1}^\infty $, if
 ${\lambda _{{n_0}}} = 0$.
\end{itemize}
are Riesz bases in ${L^2}[0, \pi]$.
\item[3.]
Let $\alpha  = \pi$, $\beta  = 0$. The answer is the same as in case 2.
\end{itemize}
\end{theorem}

The Riesz basicity of the systems of functions of sines and cosines
in ${L^2}[0, \pi]$ has been studied in many papers
(see, for example,
\cite{Freiling_Yurko:2001,He_Volkmer:2001, Levinson:1940,Moiseev:1984,
Moiseev:1987,Sedletskii:1988}) and is also associated with
Riesz basicity in ${L^2}[- \pi, \pi]$  the systems of the form
 $\{ {{e^{i{\lambda _n}x}}}\}_{n =  - \infty }^\infty $
 (see, e.g. \cite{Kadec:1964,Levin:1956,Paley_Wiener:1934}).
Completeness and Riesz basicity of systems of sines and cosines are
used in many related areas of mathematics, in particular,
in solutions of the inverse problems in spectral theory of operators
(see, e.g.
\cite{Freiling_Yurko:2001,Gasymov_Levitan:1964,Gelfand_Levitan:1951,
Levitan_Sargsyan:1988,Marchenko:1952}).

This article is organized as follows.
In section~\ref{sec2} we give some necessary information and the results
of \cite{He_Volkmer:2001}, which are more similar to ours.
In section~\ref{sec3} we prove theorems \ref{thm1.1} and \ref{thm1.2}.

\section{Preliminaries}\label{sec2}

\subsection*{Eigenvalues of the problem $L(q,\alpha, \beta )$}
The dependence of the eigenvalues of the Sturm-Liouville problem on parameters
$\alpha$ and $\beta$ from the boundary conditions \eqref{eq1.2}
and \eqref{eq1.3} was investigated in \cite{Harutyunyan:2008},
 where the following theorem was proved.

 \begin{theorem} \label{thm2.1}
  The smallest eigenvalue has the property
 \begin{equation} \label{eq2.1}
 \lim_{\alpha  \to 0} {\mu _0}(q, \alpha, \beta )  = - \infty,\quad
 \lim_{\beta  \to \pi } {\mu _0}(q, \alpha, \beta) = - \infty.
 \end{equation}
For  eigenvalues ${\mu _n}(q, \alpha, \beta)$, $n \ge 2$, the formula
 \begin{equation} \label{eq2.2}
{\mu _n}(q, \alpha, \beta) = {[n + {\delta _n}(\alpha ,\beta )]^2} + [q] + {r_n}(q, \alpha, \beta)
 \end{equation}
holds, where $[q] = \frac{1}{\pi }\int_0^\pi {q(x)dx}$,
\begin{equation} \label{eq2.3}
\begin{aligned}
 \delta _n(\alpha ,\beta )
 &= \frac{1}{\pi }\Big[\arccos \frac{\cos \alpha }
 {\sqrt {[n + {\delta _n}(\alpha ,\beta )]^2 \sin^2\alpha
 +\cos^2\alpha}} \\
 &\quad  - \arccos \frac{\cos \beta }
 {\sqrt {[n + {\delta _n}(\alpha ,\beta )]^2 \sin ^2 \beta
 + \cos ^2\beta }} \Big]
 \end{aligned}
 \end{equation}
and ${r_n}(q, \alpha, \beta ) = o(1)$, when $n \to \infty $, uniformly in
$\alpha, \beta \in [0, \pi]$ and $q$ from the bounded subsets of
$L_R^1[0, \pi]$ (we will write $q \in BL_R^1[0, \pi]$).
\end{theorem}

Note that the formula \eqref{eq2.2} is the generalization of the
asymptotic formulas known prior to \cite{Harutyunyan:2008}
for the eigenvalues of the Sturm-Liouville problem
(see \cite{Freiling_Yurko:2001,Levitan_Sargsyan:1988,Marchenko:1952,
Marchenko:1977}. More detailed table of the asymptotic formulas for
eigenvalues of the problem $L(q, \alpha, \beta)$ is in \cite{Marchenko:1952}).
From \eqref{eq2.2} for ${\lambda _n}(q, \alpha, \beta)$
$({\mu _n} = \lambda _n^2)$ we obtain the formula
 \begin{equation} \label{eq2.4}
 {\lambda _n}(q, \alpha, \beta) = n + {\delta _n}(\alpha, \beta)
+ \frac{{[q]}}{{2[n + {\delta _n}(\alpha, \beta)]}}
+ {l_n}(q, \alpha, \beta)
 \end{equation}
where ${l_n} = {l_n}(q, \alpha, \beta) = o({n^{ - 1}})$ when
$n \to \infty$ uniformly for all $q \in BL_R^1[0,\pi ]$ and
$\alpha, \beta \in [0,\pi]$. From \eqref{eq2.3} easily follows that
${\delta _n}(\alpha, \beta) = O({n^{ - 1}})$ for $\alpha, \beta \in (0,\pi)$;
${\delta _n}(\alpha, \beta) = \frac{1}{2} + O({n^{- 1}})$ for
$\alpha = \pi, \beta \in (0, \pi)$ and $\alpha \in (0, \pi), \beta = 0$; and
${\delta _n}(\pi, 0) = 1$ for all $n=2,3,\dots$. Thus, we distinguish 3 cases:
\begin{itemize}
\item[1.]
$\alpha, \beta \in (0, \pi)$; i.e. the interior points of the square
$[0, \pi] \times [0, \pi]$, where
${\lambda _n} = {\lambda _n}(q, \alpha, \beta)$ have the asymptotic property
${\lambda _n} = n + O({n^{ - 1}})$,

\item[2.]
$\alpha  = \pi$, $\beta \in (0, \pi)$ or $\alpha \in (0, \pi)$, $\beta = 0$
(i.e. right and bottom edges of the square $[0, \pi] \times [0, \pi]$),
where ${\lambda _n}$ have the asymptotic property
 ${\lambda _n} = n + \frac{1}{2} + O({n^{ - 1}})$,

\item[3.]
$\alpha  = \pi, ~\beta  = 0$, where
${\lambda _n}(q,\pi ,0) = n + 1 + O({n^{ - 1}})$.
\end{itemize}

\subsection*{Riesz bases}
The following three definitions  and two lemmas are taken from
\cite{Freiling_Yurko:2001}.
Equivalent definitions and statements are available in other studies
(see, e.g. \cite{Gohberg_Krein:1965,He_Volkmer:2001,Levin:1956}).

\begin{definition} \label{def2.1} \rm
A basis $\{ {f_j}\} _{j = 1}^\infty$  of a separable Hilbert space $H$
is called a Riesz basis if it is derived from an orthonormal basis
$\{ {e_j}\} _{j = 1}^\infty$  by linear bounded invertible operator
$A$, i.e., if ${f_j} = A{e_j}$, $j=1,2,\dots$.
\end{definition}

\begin{definition} \label{def2.2} \rm
Two sequences of elements  $\{ {f_j}\} _{j = 1}^\infty$ and
$\{ {g_j}\} _{j = 1}^\infty$ from $H$ are called quadratically close
if $\sum_{j = 1}^\infty  {{{\|{f_j} - {g_j}\|}^2}}  < \infty $.
\end{definition}

\begin{definition} \label{def2.3} \rm
 A sequence $\{ {g_n}\} _{n = 0}^\infty$ is called $\omega$-linearly
independent, if the equality $\sum_{n = 0}^\infty  {{c_n}{g_n}}  = 0$
is possible only when ${c_n} = 0$ for $n = 0,1,2,\dots$.
\end{definition}

\begin{lemma} \label{lem2.1}
Let $\{ {f_n}\} _{n = 0}^\infty $ be a Riesz basis in $H$,
$\{ {f_n}\} _{n = 0}^\infty $ and $\{ {g_n}\} _{n = 0}^\infty $
are quadratically close. If $\{ {g_n}\} _{n = 0}^\infty $ is $\omega$-linearly
independent, then $\{ {g_n}\} _{n = 0}^\infty $ is a Riesz basis in $H$.
\end{lemma}

\begin{lemma} \label{lem2.2}
Let $\{ {f_n}\} _{n = 0}^\infty $ be a Riesz basis in $H$,
$\{ {f_n}\} _{n = 0}^\infty $ and $\{ {g_n}\} _{n = 0}^\infty $
are quadratically close. If $\{ {g_n}\} _{n = 0}^\infty $ is complete in $H$,
then $\{ {g_n}\} _{n = 0}^\infty $ is $\omega$-linearly independent
(and therefore, is a Riesz basis in $H$).
\end{lemma}

The following two theorems are proved in \cite{He_Volkmer:2001}.

\begin{theorem} \label{thm2.2}
 Let $\{ {\lambda _n}\} _{n = 0}^\infty $ be a sequence of nonnegative
numbers with the property that ${\lambda _k} \ne {\lambda _m}$ for
$k \ne m$  and of the form ${\lambda _n} = n + \delta  + {\delta _n}$,
with ${\delta _n} \in [-l,l]$ for sufficiently large $n$, where
the constants $\delta  \in [0,\frac{1}{2}]$ and  $l \in (0,\frac{1}{4})$
satisfy ${[1 + \sin (2\pi \delta )]^{\frac{1}{2}}}(1 - \cos (\pi l))
+ \sin (\pi l) < 1$. Then $\{ \cos(\lambda _n x)\} _{n = 0}^\infty $
is a Riesz basis in ${L^2}[0, \pi]$.
\end{theorem}

\begin{theorem} \label{thm2.3}
Let $\{ {\lambda _n}\} _{n = 1}^\infty$ be a sequence of positive numbers
of the form ${\lambda _n} = n - \delta  + {\delta _n}$, having the same
properties as in Theorem \ref{thm2.2}.
Then $\{ \sin(\lambda_n x)\} _{n = 1}^\infty $is a Riesz basis in ${L^2}[0,\pi ]$.
\end{theorem}

It follows from  \eqref{eq2.4} that the only circumstance, (essentially)
preventing us to apply theorems \ref{thm2.2} and \ref{thm2.3}
for proving Riesz basicity of systems
$\{ \cos(\lambda _n(q, \alpha, \beta)x)\} _{n = 0}^\infty $ and
$\{ \sin(\lambda _n(q,\alpha,\beta)x)\} _{n = 0}^\infty$,
is that among the eigenvalues ${\mu _n} = \lambda _n^2(q, \alpha, \beta)$
may be negative (see \eqref{eq2.1}), and accordingly among
${\lambda _n}(q, \alpha, \beta)$ may be (in a finite number) pure imaginary ones.
 Can these $\lambda_n$ interfere the Riesz basicity of the mentioned systems?
 Our answer is contained in theorems \ref{thm1.1} and \ref{thm1.2}.

\section{Proofs of main Theorems} \label{sec3}

We will start with a lemma, which is an analogue of
\cite[Lemma 4]{He_Volkmer:2001}.

\begin{lemma} \label{lem3.1}
Let $\{ {\nu _n^2}\} _{n = 0}^\infty $ and
 $\{ {\lambda _n^2}\} _{n = 0}^\infty$ be two real sequences such that
${\nu _k^2} \ne {\nu _m^2}$ and ${\lambda _k^2} \ne {\lambda _m^2}$,
for $k \ne m$, and among which only a finite number of members
$({\nu _0^2}, {\nu _1^2},\dots ,{\nu _{{n_1}}^2}; {\lambda _0^2},
{\lambda _1^2},\dots, {\lambda _{{n_2}}^2})$ can be negative,
and the sequences are enumerated in increasing order
$({\nu _0^2} < {\nu _1^2} < \dots  < {\nu _n^2} < \dots ;
{\lambda _0^2} < {\lambda _1^2} < \dots  < {\lambda _n^2} < \dots )$.
Let $\{ {\nu _n}\} $ and $\{ {\lambda _n}\} $ have the asymptotic properties
\begin{gather} \label{eq3.1}
{\nu _n} = n + \delta  + O({n^{ - 1}}), \quad 0 \le \delta  \le 1,\\
 \label{eq3.2}
{\lambda _n} = n + {\delta _n}(\alpha, \beta) + O({n^{ - 1}}),
\end{gather}
when $n \to \infty$ and, furthermore,
\begin{equation} \label{eq3.3}
\sum_{n = 0}^\infty  {{|\lambda _n - \nu _n |^2}} < \infty.
\end{equation}
Then $\{ \cos(\nu_n x)\} _{n = 0}^\infty $ is a Riesz basis in
${L^2}[0, \pi]$ if and only if $\{ \cos(\lambda_nx)\} _{n = 0}^\infty $
is a Riesz basis in ${L^2}[0, \pi]$.
\end{lemma}

\begin{proof}
Set ${f_n}(x) = \cos(\nu_n x)$ and ${g_n}(x) = \cos(\lambda_nx)$,
$n=0,1,2\dots$. Assume $\{ {f_n}\} _{n = 0}^\infty$ is a Riesz basis in
 ${L^2}[0, \pi]$. Since for real numbers ${\nu _n}$ and $\lambda _n$,
\begin{align*}
|\cos(\nu_n x) - \cos(\lambda_nx)|
&=|2\sin \frac{{({\lambda _n} - {\nu _n})x}}{2}\sin \frac{{({\lambda _n}
+ {\nu _n})x}}{2}| \\
&\le 2 | {\sin \frac{{({\lambda _n} - {\nu _n})x}}{2}} |
\le | {{\nu _n} - {\lambda _n}} |x \le \pi | {{\nu _n} - {\lambda _n}} |,
\end{align*}
 we obtain that
\[
{\| {\cos(\nu_n x) - \cos(\lambda_nx)} \|^2}
= \int_0^\pi  {{| \cos(\nu_n x) - \cos(\lambda_nx)| ^2}} dx
\le {\pi ^3}{| {{\nu _n} - {\lambda _n}}|^2}.
\]
 Therefore,
\begin{align*}
\sum_{n = 0}^\infty  {{\| {{f_n} - {g_n}}\|}^2}
&= \sum_{n = 0}^{{n_0}} {{{\| {{f_n} - {g_n}}\|}^2}}
+ \sum_{n = {n_0} + 1}^\infty  {{{\| {{f_n} - {g_n}} \|}^2}} \\
&\le {M_0}+{\pi ^3}\sum_{n = {n_0} + 1}^\infty  {{| {{\lambda _n}
- {\nu _n}} |}^2} < \infty ;
\end{align*}
 i.e., $\{ {f_n}\} _{n = 0}^\infty $ and $\{ {g_n}\} _{n = 0}^\infty $
are quadratically close $({n_0} = \max \{ {n_1},{n_2}\} )$.
 According to Lemma \ref{lem2.1}, to prove the Riesz basicity of the
system $\{ {g_n}\} _{n = 0}^\infty $ it is enough to prove its $\omega$-linearly
independence. Assume the contrary, i.e. let there is a sequence
$\{ {c_n}\} _{n = 0}^\infty  \in {l^2}$, not identically zero, such that
\begin{equation} \label{eq3.4}
\sum_{n = 0}^\infty  {{c_n}} {g_n} = 0.
\end{equation}
Let $\lambda  \in\mathbb{C}$ be such that
$\lambda  \ne  \pm {\lambda _n}$, $n=0,1,2,\dots$, and define the function
\begin{equation} \label{eq3.5}
g(x) = \sum_{n = 0}^\infty  {\frac{{{c_n}}}{{\lambda _n^2 - {\lambda ^2}}}}
 {g_n}(x).
\end{equation}
It follows from  \eqref{eq2.2} that this series is uniformly convergent
for $x \in [0, \pi]$. Similarly, the series
\begin{equation} \label{eq3.6}
g'(x) =  - \sum_{n = 0}^\infty
{\frac{{{c_n}{\lambda _n}}}{{\lambda _n^2 - {\lambda ^2}}}} \sin(\lambda_n x)
\end{equation}
converges uniformly on $[0,\pi ]$. Since ${g''_n} =  - \lambda _n^2{g_n}$,
 we have (note, that here we repeat the proof of \cite{He_Volkmer:2001})
\[
\sum_{n = 0}^m {\frac{{{c_n}}}{{\lambda _n^2 - {\lambda ^2}}}} {g''_n}(x)
= - \sum_{n = 0}^m {\frac{{{c_n}\lambda _n^2}}{{\lambda _n^2
 - {\lambda ^2}}}} {g_n}(x) = - \sum_{n = 0}^m {{c_n}{g_n}} (x)
  - {\lambda ^2}\sum_{n = 0}^m {\frac{{{c_n}}}{{\lambda _n^2
 - {\lambda ^2}}}} {g_n}(x).
\]
Taking into account \eqref{eq3.4}, we conclude that the sequence on
the left-side of the last equality converges in ${L^2}[0,\pi ]$ to
 $- {\lambda ^2}g(x)$, when $m \to \infty $. This implies that $g$
is twice differentiable and satisfies the differential equation
$- g''(x) = {\lambda ^2}g(x),x \in (0,\pi )$, and initial conditions
(see \eqref{eq3.5} and \eqref{eq3.6}):
\begin{equation} \label{eq3.7}
g(0) = h(\lambda ) = \sum_{n = 0}^\infty  {\frac{{{c_n}}}{{\lambda _n^2
- {\lambda ^2}}}}, ~g'(0) = 0;
\end{equation}
i.e., $g$ is the solution of the corresponding Cauchy problem,
which is unique and given by the formula
\begin{equation} \label{eq3.8}
g(x) = h(\lambda )\cos(\lambda x).
\end{equation}
The function $h(\lambda )$ defined by  \eqref{eq3.7} is meromorphic,
and taking into account that
 $\{ {c_n}\} _{n = 0}^\infty  \ne \{ 0\} _{n = 0}^\infty$,
is not an identically zero function. Then it has no more than countable
 number of isolated zeros.
If $h(\lambda ) \ne 0$, \eqref{eq3.5} and \eqref{eq3.8} show that
$\cos(\lambda x)$ belongs to the closed linear span of the system
$\{ {g_n}\} _{n = 0}^\infty $ in ${L^2}[0, \pi]$. Since $\cos(\lambda x)$
is a continuous function of $(\lambda ,x)$, we obtain that $\cos(\lambda x)$
belongs to closed linear span of the system $\{ {g_n}\} _{n = 0}^\infty$
for all $\lambda \in\mathbb{C}$. Particularly, the all $\cos(nx)$,
$n = 0,1,2,\dots$ belong to the closed linear span of the
system $\{ {g_n}\} _{n = 0}^\infty$, so the system $\{ {g_n}\} _{n = 0}^\infty $
is a complete system in ${L^2}[0, \pi]$. From Lemma \ref{lem2.2} follows
the $\omega$-linearly independence of the system $\{ {g_n}\} _{n = 0}^\infty$;
i.e. we come to contradiction, and the Riesz basicity of the system
$\{ {g_n}\} _{n = 0}^\infty $ is proved. If we assume, that
$\{ {g_n}\} _{n = 0}^\infty$ is a Riesz basis, then similarly we can prove
the Riesz basicity of the system $\{ {f_n}\} _{n = 0}^\infty$.
Lemma \ref{lem3.1} is proved.
\end{proof}

Let us now turn to the proof of the theorem \ref{thm1.1}.
 Let us start from the first case:
 $\alpha, \beta \in (0, \pi)$.
Let us take in Lemma \ref{lem3.1} ${\nu _n} = n, {f_n}(x) = \cos(nx)$, and
${g_n}(x) = \cos(\lambda _n(q, \alpha, \beta)x)$, $n=0,1,2,\dots$.
In this case ${\lambda _n}(q, \alpha, \beta) = n + O({n^{ - 1}})$, and,
therefore,  \eqref{eq3.3} holds; i.e.,
 $\{ {f_n}\} $ and $\{ {g_n}\} $ are quadratically close.
Since $\{ {f_n}\} _{n = 0}^\infty  = \{ \cos(nx)\} _{n = 0}^\infty$
is a Riesz basis, then from the Lemma \ref{lem3.1} follows the Riesz
basicity of the system $\{ \cos(\lambda _n(q, \alpha, \beta)x)\} _{n = 0}^\infty$.

In the second case in Lemma \ref{lem3.1} we take
${\nu _n} = n + \frac{1}{2}$, ${f_n}(x) = \cos ((n + \frac{1}{2})x)$, and
${g_n}(x) = \cos(\lambda _n(q, \alpha, \beta )x)$, $n=0,1,2,\dots$.
In the second case ${\lambda _n}(q, \alpha, \beta )
= n + \frac{1}{2} + O({n^{ - 1}})$ and therefore again holds \eqref{eq3.3};
i.e. quadratically closeness. As
$\{ \cos ((n + \frac{1}{2})x)\} _{n = 0}^\infty $ is the system of eigenfunctions
of the Sturm-Liouville problem $L(0,\frac{\pi }{2},0)$, it is an orthogonal
basis in ${L^2}[0, \pi]$ (and, particularly, is a Riesz basis).
From the Lemma \ref{lem3.1} follows the Riesz basicity of the system
$\{ \cos(\lambda _n(q, \alpha, \beta)x)\} _{n = 0}^\infty$ in this case.

In the third case in Lemma \ref{lem3.1} we take
${\nu _n} = n + 1$ $(\delta  = 1)$,
${f_n}(x) = \cos ((n + 1)x)$ and
${g_n}(x) = \cos ({\lambda _n}(q, \pi, 0)x)$, $n=0,1,2,\dots$.
If we assume that $\{ {g_n}\} _{n = 0}^\infty $ is a Riesz basis,
then from the asymptotic property
${\lambda _n}(q, \pi, 0) = n + 1 + O({n^{ - 1}})$ and Lemma \ref{lem3.1}
follows the Riesz basicity of the system
$\{ {f_n}\} _{n = 0}^\infty  = \{ \cos ((n + 1)x)\} _{n = 0}^\infty$,
which is incorrect, since it is even not complete.
Therefore, $\{ \cos {\lambda _n}(q,\pi ,0)x\} _{n = 0}^\infty $
does not form a Riesz basis. But adding to this system a function
$f(x)=\cos(\lambda x)$, where $\lambda ^2 \neq \lambda _n^2$ for every
$n=0,1,2,\dots$ and noticing that the system
$\{ f(x)\}  \cup \{ \cos(\lambda _n(q, \pi, 0)x)\} _{n = 0}^\infty$
is $\omega$-linearly independent and quadratically close to the system
$\{ \cos(nx)\} _{n = 0}^\infty$, according to the Lemma \ref{lem3.1},
 we get its Riesz basicity. Theorem \ref{thm1.1} is proved.

\begin{lemma} \label{lem3.2}
Let $\{ {\nu _n^2}\} _{n = 0}^\infty$ and $\{ {\lambda _n^2}\} _{n = 0}^\infty $
are the same as in Lemma \ref{lem3.1}. Then
$\{ \sin(\nu _n x)\} _{n = 0}^\infty$ is a Riesz basis in ${L^2}[0, \pi]$
if and only if $\{ \sin(\lambda_n x)\} _{n = 0}^\infty$ is a Riesz basis
 in ${L^2}[0, \pi]$.
\end{lemma}

\begin{proof}
Set ${f_n}(x) = \sin(\nu _n x)$, ${g_n}(x) = \sin(\lambda_n x)$, $n=0,1,2,\dots$.
 Quadratic closeness of the systems $\{ {f_n}\} _{n = 0}^\infty$ and
$\{ {g_n}\} _{n = 0}^\infty$ can be showed in the same way as in
Lemma \ref{lem3.1}. Function $g(x)$ (see \eqref{eq3.5}) in this case
is the solution of the Cauchy problem $ - g'' = {\lambda ^2}g$,
$g(0) = 0$,
$g'(0) = {h_1}(\lambda ) = \sum_{n = 0}^\infty
{\frac{{{c_n}{\lambda _n}}}{{\lambda _n^2 - {\lambda ^2}}}}$, and,
therefore, has the form
$g(x) = {h_1}(\lambda ) \sin(\lambda x)/\lambda $.
From the continuity of $\sin(\lambda x)/\lambda$
as a function of two variables $(\lambda, x)$ follows that the equality
\begin{equation} \label{eq3.9}
\frac{{\sin (\lambda x)}}{\lambda }
= \frac{1}{{{h_1}(\lambda )}}\sum_{n = 0}^\infty
{\frac{{{c_n}}}{{\lambda _n^2 - {\lambda ^2}}}} \sin(\lambda_n x)
\end{equation}
holds not only when ${h_1}(\lambda ) \ne 0$, but also for all
$\lambda  \in\mathbb{C}$. Hence \eqref{eq3.9} is right for
$\lambda = 1,2,3,\dots$; i.e., the all elements of the orthogonal basis
$\{ \sin(nx)\} _{n = 1}^\infty$ are in the closed linear span of the
system $\{ \sin(\lambda_n x)\} _{n = 0}^\infty$; i.e., the system
 $\{ \sin(\lambda_n x)\} _{n = 0}^\infty$ is complete in ${L^2}[0, \pi]$.
The rest of the proof is as in Lemma \ref{lem3.1}.
\end{proof}

Now the proof of Theorem \ref{thm1.2} is: In stated in following cases:
\begin{itemize}
\item[(1.a)] we take ${\nu _n} = n + 1$ and accordingly
$\{ {f_n}(x)\} _{n = 0}^\infty  = \{ \sin ((n + 1)x)\} _{n = 0}^\infty
= \{ \sin(nx)\} _{n = 1}^\infty$ and
$\{ {g_n}(x)\} _{n = 1}^\infty  = \{ \sin(\lambda_n x)\} _{n = 1}^\infty$,
as stated in Theorem \ref{thm1.2}.
Since $\{ {f_n}\} _{n = 0}^\infty $ is a Riesz basis
(and even an orthogonal basis) and from the asymptotic property
 ${\lambda _n} = n + O({n^{ - 1}})$ it follows that
$\{ {f_n}\}$ and $\{ {g_n}\}$ are quadratically close,
therefore from Lemma \ref{lem3.2} follows the Riesz basicity of the
system $\{ \sin(\lambda_n x)\} _{n = 1}^\infty$.

\item[(1.b)] also $\{ {f_n}(x)\} _{n = 1}^\infty
= \{ \sin(nx)\} _{n = 1}^\infty$ and the system
\[
\{ \sin(\lambda_n x)\} _{n = 0}^{{n_0} - 1} \cup
\{ \sin(\lambda_n x)\} _{n = {n_0} + 1}^\infty
\]
is again quadratically close to $\{ {f_n}\} _{n = 1}^\infty$.

\item[(2.a)] we take ${\nu _n} = n + \frac{1}{2}$, accordingly,
${f_n}(x) = \sin ((n + \frac{1}{2})x)$, and
${g_n}(x) = \sin(\lambda_n x)$, $n=0,1,\dots$. Since
$\{ \sin ((n + \frac{1}{2})x)\} _{n = 0}^\infty $ is the system of
eigenfunctions of the self-adjoint problem $L(0, \pi, \frac{\pi}{2})$,
than it is an orthogonal basis in ${L^2}[0, \pi]$. The asymptotic property
 ${\lambda _n} = n + \frac{1}{2} + O({n^{ - 1}})$ ensures the quadratically
closeness of the systems $\{ {f_n}\} _{n = 0}^\infty$ and
$\{ {g_n}\} _{n = 0}^\infty$, therefore in this case the Riesz basicity
of the system $\{ {g_n}\} _{n = 0}^\infty$ is proved.

\item[(2.b)] again ${f_n}(x) = \sin ((n + \frac{1}{2})x)$,
$n=0,1,\dots$, and $\{ {g_n}\} $ is different from the case (2.a)
 with only one element ${g_{{n_0}}}$, which has not any effect on
quadratically closeness of the systems $\{ {f_n}\} _{n = 0}^\infty$ and
$\{ {g_n}\} _{n = 0}^\infty$.

\item[(3)] we take ${\nu _n} = n + 1$ and ${f_n}(x) = \sin ((n + 1)x)$,
$n=0,1,2,\dots$; i.e., $\{ {f_n}(x)\} _{n = 0}^\infty
= \{ \sin(nx)\} _{n = 1}^\infty$. The rest is followed from the asymptotic
property ${\lambda _n}(q, \pi, 0) = n + 1 + O({n^{ - 1}})$, if we take
${g_n}(x) = \sin ({\lambda _n}(q, \pi, 0)x)$,  $n=0,1,\dots$.
\end{itemize}
Therefore, theorem \ref{thm1.2} is proved.

\begin{remark} \rm
From lemmas \ref{lem3.1} and \ref{lem3.2} it easily follows
that $\{ \cos(\lambda _n (q, \alpha, \beta)x)\} _{n = 0}^\infty$
is a Riesz basis in ${L^2}[0,\pi ]$ if and only if
 $\{ \cos (\lambda _n(0, \alpha, \beta)x)\} _{n = 0}^\infty $
 is a Riesz basis in ${L^2}[0, \pi]$. Similarly for sines.
This means that the stability of Riesz basicity is not affected
by adding the potential $q(\cdot)$.
\end{remark}

\subsection*{Acknowledgments}
The authors are deeply indebted to the anonymous referee for the valuable
suggestions and comments which improved this manuscript.

\end{document}